\def\IH{{\Bbb H}}

\def\IS{{\Bbb S}} 
 
\def\IK{{\Bbb K}}

\def\IC{\Bbb C} 
 
\def\ID{{\Bbb D}}

\def\zbar{{\overline{z}}} 
 
\def\wbar{{\overline{w}}}

\documentclass[12pt]{article} 

\usepackage[psamsfonts]{amssymb} 
\usepackage{amsfonts,amsmath} 
\usepackage{tikz-cd}

\usepackage{epsfig,multicol}

\newtheorem{theorem}{Theorem}
\newtheorem{lemma}{Lemma}
\newtheorem{corollary}{Corollary}

\title{On the uniqueness of extremal mappings of finite distortion.}
\author{Gaven Martin \& Cong Yao \thanks{
Work of both authors partially supported by the New Zealand Marsden Fund.
\newline
Institute for Advanced Study, 
Massey University,  Auckland,
New Zealand.
\newline
email: g.j.martin@massey.ac.nz \&
c.yao@massey.ac.nz
\newline
{\bf Keywords.} Quasiconformal,  finite distortion,  uniqueness, extremal mappings, calculus of variations
\newline
{\bf MSC Subject.}  30C62 31A05 49J10  } }
\date{}
\begin{document}
\maketitle
\begin{abstract} For an arbitrary convex function $\Psi:[1,\infty)  \to [1,\infty)$, we consider uniqueness in the following two related extremal problems:  

 \noindent{\bf Problem A}  {\em (boundary value problem):} Establish the existence of,  and describe the mapping $f$, achieving
\[ 
\inf_f \Big\{ \int_\ID \Psi(\IK(z,f))\; dz : f:\bar\ID \to \bar\ID \; \mbox{a homeomorphism in $W^{1,1}_{0}(\ID)+f_0$}  \Big\}.  
\] 
Here the data $f_0:\bar\ID \to \bar\ID$ is a homeomorphism of finite distortion with $\int_\ID \Psi(\IK(z,f_0))\; dz<\infty$ -- a barrier.  
Next,  given two homeomorphic Riemann surfaces $R$ and $S$ and data $f_0:R \to S$ a diffeomorphism.

 \noindent{\bf Problem B}  {\em (extremal in homotopy class):} Establish the existence of,  and describe the mapping $f$, achieving
\[ \inf_f \Big\{ \int_R \Psi(\IK(z,f))\; \;d\sigma(z) : \mbox{$f$ a homeomorphism homotopic to $f_0$}  \Big\}.  \]  

There are two basic obstructions to existence and regularity.  These are first, the existence of an Ahlfors-Hopf differential and second that the minimiser is a homeomorphism.  When these restrictions are met (as they often can be) we show uniqueness is assured. These results are established through a generalisation the classical Reich-Strebel inequalities to this variational setting.
\end{abstract}

\section{Introduction.}  There are two usual approaches to the Teichm\"uller theory of Riemann surfaces, the classical approach through extremal quasiconformal mappings of Teichm\"uller and Ahlfors \cite{Ahlfors} and the approach through harmonic mappings initiated by   Fischer\&Tromba, Jost Wolf and Wolpert \cite{FT,Jost,Wolf,Wolf2},  the latter approach building on important earlier work of many people.  Both approaches have their appeal and many wonderful results are proved using either or both approaches -- for instance the solution of the Nielsen realisation problem \cite{Kerckhoff,Wolpert}.

\medskip

Our recent work \cite{MY1,MY2,MY3} seeks to unite these different approaches through the calculus of variations and intermediate Teichm\"uller theories, in particular the $L^p$ and  $exp_q$ Teichm\"uller theories, where $\Psi(t)\approx t^p$, $1\leq p<\infty$ and $\Psi(t)\approx e^{t,q}$, $0<q<\infty$ in Problems A \& B below.  In both cases as $p\to\infty$ we obtain an extremal quasiconformal mapping as the limits of minimisers (despite the fact this sequence may not consist of homeomorphism in the $L^p$ case).  While as $p\to 1$ or $q\to 0$ the minimisers converge to a harmonic mapping. Our results build on earlier work concerning extremal mappings of finite distortion \cite{AIMO,AIM,IMO,MY4}.  

\section{The problems.}

The distortion function we will minimise is 
\[ \IK(z,f) = \frac{|f_z|^2+|f_\zbar|^2}{|f_z|^2-|f_\zbar|^2} =\frac{\|Df\|^2}{J(z,f)} \]
It was already realised by Ahlfors \cite{Ahlfors} in his rigorous proof of Teichm\"uller's Theorem that this is the correct distortion to use as opposed the the distortion $K(z,f)=\frac{|f_z|+|f_\zbar|}{|f_z|-|f_\zbar|} =\frac{|Df|^2}{J(z,f)}$.  The functional relation $\IK=\frac{1}{2}(K+1/K)$ shows the $L^{\infty}$ minimisers to be the same.

\medskip

For an arbitrary convex  increasing  function $\Psi:[1,\infty)  \to [1,\infty)$ with $\Psi(t)\geq t$, we consider uniqueness in the following two related extremal problems:  

\medskip

 \noindent{\bf Problem A}  {\em (boundary value problem):} Find and describe the mapping $f$ achieving
\begin{equation} 
\inf_f \Big\{ \int_\ID \Psi(\IK(z,f))\; \lambda(z) \;dz : f:\bar\ID \to \bar\ID \; \mbox{a homeomorphism in $W^{1,1}_{0}(\ID)+f_0$}  \Big\}.  
\end{equation}  
Here the data $f_0:\bar\ID \to \bar\ID$ is a homeomorphism of finite distortion with $\int_\ID \Psi(\IK(z,f_0))\,\lambda(z)\; dz<\infty$  and $\lambda\geq 1$ a smooth weight.  

\medskip

Next,  given two homeomorphic Riemann surfaces $R$ and $S$ and data $f_0:R \to S$ a diffeomorphism.

\medskip

 \noindent{\bf Problem B}  {\em (extremal in homotopy class):} Find and describe the mapping $f$ achieving
\[ \inf_f \Big\{ \int_R \Psi(\IK(z,f))\; \;d\sigma(z) : \mbox{$f$ a homeomorphism homotopic to $f_0$}  \Big\}.  \]  

There are two basic obstructions to existence and regularity \cite{MY1,MY2}.  These are first, the existence of an Ahlfors-Hopf differential - a holomorphic quadratic differential arising from certain variational equations, and second that the minimiser is a homeomorphism.  However both of these can be achieved for locally quasiconformal mappings.

\begin{theorem}\label{main1} A locally quasiconformal minimiser to Problem A or Problem B is a diffeomorphism and is unique among homeomorphic minimisers.
\end{theorem}

\noindent{\bf Remark 1.} In Theorem \ref{main1} one can make rather more precise statements depending on the growth of $\Psi$. For instance
\begin{itemize}
\item If $\Psi(t)\approx t^p$, $p\geq 1$, and $\IK(z,f) \in L^r_{loc}$ for any $r>1+p$ (notice apriori $\IK(z,f) \in L^p$), then the minimiser to Problem A or Problem B is a diffeomorphism and is unique among homeomorphic minimisers. 
\item If $\Psi(t)\approx e^{pt}$, $p>0$, and $e^{r\IK(z,f)} \in L^1_{loc}$ for any $r>p$ (notice apriori $e^{p\IK(z,f)} \in L^1$), then the minimiser to Problem A is a diffeomorphism and is unique. 
\end{itemize}

Problem B admits a complete solution if $\Psi$ grows fast enough.

\begin{theorem}
If $\Psi(t)\geq e^{\epsilon\,  t}$, for some $\epsilon>0$, then a minimiser to Problem B exists, is a diffeomorphism and is unique. 
\end{theorem}

In both cases as $p\to\infty$ we recover  Ahlfors' approach to Teichm\"uller theory in the sense that for fixed data the uniform limit of the $p$-extremal mappings is a Teichm\"uller extremal,  and as $p\to 1$ (for the $L^p$ case) or as $p\to 0$ (for the $exp_q$ case) we recover the extremal harmonic mapping.  In either case,  as $p$   varies we get smoothly varying families of degenerate elliptic Beltrami systems controlling the properties of the extremal mappings and extremal mappings are harmonic in a metric induced by their distortion.  

Our results reveal the following very curious dichotomy - analogous to the issues for extremal quasiconformal mappings (see \cite{Hamilton, RS2,RS3, BLMM,Yao}).  Given fixed data -- boundary values or a homotopy class and so forth -- we find in 
\begin{itemize}
\item  $L^p$-case:  There is always an Alhfors-Hopf quadratic differential, but minimisers might not be homeomorphisms.
\item $exp$-case:  There is always a homeomorphic minimiser,  but there might not be an Alhfors-Hopf quadratic differential.
\end{itemize}

\section{Finite distortion functions.}
Let $\Omega$ be a planar domain. A mapping $f:\Omega\to\IC$ has finite distortion if 
\begin{enumerate}
\item $f\in W^{1,1}_{loc}(\Omega)$,  the Sobolev space of functions with locally integrable first derivatives,
\item the Jacobian determinant $J(z,f)\in L^{1}_{loc}(\Omega)$, and 
\item there is a measurable function $ \IK(z,f)\geq 1$, finite almost everywhere, such that 
 \begin{equation}\label{1.1}
 \| Df(z) \|^2 \leq \IK(z,f) \, J(z,f), \hskip10pt \mbox{ almost everywhere in $\Omega$}.
 \end{equation}
\end{enumerate}
The smallest such function $\IK(z,f)$ is called the distortion of $f$. One obtains the class of quasiconformal functions if $f$ is a homeomorphism and $\|\IK(z,f)\|_{L^\infty(\Omega)}<\infty$,  in which case $f\in W^{1,2}_{loc}(\Omega)$ automatically.  We recommend \cite[Chapter 20]{AIM} for the basic theory of mappings of finite distortion and the associated governing equations; degenerate elliptic Beltrami systems.  We recall that the Beltrami coefficient of $f$ is
\begin{equation}
\mu_f(z) = \frac{f_\zbar(z)}{f_z(z)},
\end{equation}
and hence
\begin{equation}
\IK(z,f) = \frac{1+|\mu_f|^2}{1-|\mu_f|^2}.
\end{equation}
For mappings between Riemann surfaces we use the same terminology and as we will only be concerned with continuous mappings of finite distortion in the article and so the subtleties of Sobolev theory in this setting will not concern us.

\section{Ahlfors-Hopf quadratic differentials and variational equations.}

It was Ahlfors \cite{Ahlfors} who first realised, in a special case, that the obvious minimisation problems suggested above are best tackled by considering a related but apparently more complicated (as it involves more nonlinearity in the integrand) problem.  This is because the variational equations for a minimiser show that a certain nonlinear combination of derivatives of a minimiser is holomorphic. 

\begin{theorem} Let $f$ be a candidate for the minimisation problems A or B (so a homeomorphic mapping of finite distortion).  Let $h=f^{-1}$.  Then in each case
\begin{equation}\label{cov}
  \int_\ID \Psi(\IK(z,f))\,\lambda(z)\; dz = \int_\ID \Psi(\IK(w,h))\,\lambda(h)\,J(w,h) \; dw  
\end{equation}
or
\[   \int_R \Psi(\IK(z,f))\; \;d\sigma_R (z) = \int_S \Psi(\IK(w,h))\; \;d\sigma_S (h) \]  
\end{theorem}

For diffeomorphisms this is simply the change of variables formula.  However in the class of mappings we seek to minimise over the result is altogether rather more subtle.  The definitive result, at least in this setting,  can be found in \cite{HK}.

\medskip

Next,  let $\varphi\in C^{\infty}_{0}(\ID)$, $\|\nabla\varphi\|_\infty<1$ (or in the case of closed Riemann surfaces $\varphi\in C^{\infty}_{0}({\cal P})$ for a convex fundamental domain).  Set $g^t=z+t\varphi$.  Then $g^t:\ID\to \ID$ is a diffeomorphism with $g^t:|\partial\ID=identity$ (we can construct an automorphic $g^t$ in the obvious way in the surface case).  Then for a sufficiently regular minimiser - should it exist, it does not need to be a homeomorphism - we can find the Euler-Lagrange, or inner-variational, equations from
\begin{equation}\label{vareqn}
0 = \frac{d}{dt}\Big|_{t=0}    \int_\ID \Psi(\IK(w,h\circ g^t))\,\lambda(h\circ g^t)\,J(w,h\circ g^t) \; dw  
\end{equation} 
and similarly in the automorphic case.  The details are given in \cite{MY1,MY2} but the key issue is the integrability of $ \Psi(\IK(w,h\circ g^t))$ when $\Psi$ grows rapidly.  It is a reasonably straightforward calculation,  finally using Weyl's lemma, to deduce that (\ref{vareqn}) for all such $\varphi$ implies the following.
\begin{theorem} Under the hypotheses of Theorem \ref{main1} or the subsequent discussion, the nonlinear combination of derivatives  
\begin{equation}
\Phi_h = \Psi'(\IK(w,h)) h_w \, h_\wbar\, \lambda(h) 
\end{equation}
is holomorphic, or in the  Riemann surface   case
\begin{equation}
\Phi_h = \Psi'(\IK(w,h)) h_w \,\overline{ h_\wbar} \; d\sigma_R(h) 
\end{equation}
is a holomorphic quadratic differential.  
\end{theorem}

See the discussion below in Remark 2. for clarification of our notation used here.  In the last case here when we lift to the universal cover $\ID$ we find
\begin{equation}\label{hypAH} \Phi_{\hat{h}} = \Psi'(\IK(w,\hat{h})) \frac{\hat{h}_w \, \overline{\hat{h}_\wbar}}{(1-|\hat{h}|^2)^2}  
\end{equation}
is a holomorphic function on $\ID$ with the property that for all $\gamma\in \pi_1(S)$,
\[ \Phi_{\hat{h}}(w)=\Phi_{\hat{h}}(\gamma(w))(\gamma'(w))^2. \]
Notice also the further invariance property of $\Phi_{\hat{h}}$ which follows from the fact M\"obius transformations are hyperbolic isometries.  If $\psi:\ID\to\ID$ is M\"obius, then
\[ \Phi_{\hat{h}}(w)=\Phi_{\psi\circ \hat{h}}(w)\]
Because of Ahlfors' discovery and the connection with quadratic differentials and the Hopf equation,  we call $\Phi_h$ an Ahlfors-Hopf quadratic differential.

\medskip

At this point it is important to make the following calculation. Let $\mu_h(w)$ be the Beltrami coefficient of $h$, $\mu_h=h_\wbar/h_w$.
\begin{eqnarray*}
\Phi_h & = &  \Psi'(\IK(w,h)) h_w \, \overline{h_\wbar}\, \lambda(h)  = \mu_h(w) \Psi'(\IK(w,h)) |h_w|^2  \lambda(h) \\
& = &  \frac{\mu_h(w)}{1-|\mu|^2}\,  \Psi'(\IK(w,h)) |h_w|^2(1-|\mu|^2)  \lambda(h) \\
& = &  \frac{\mu_h(w)}{1+|\mu|^2}\,  \IK(w,h) \Psi'(\IK(w,h)) J(w,h)  \lambda(h).
\end{eqnarray*}
The term $\frac{\mu_h(w)}{1+|\mu|^2}$ is never more than a half. The key observation is that we have apriori from (\ref{cov}) that
\[  \Psi(\IK(w,h)) J(w,h) \lambda(h) \in L^1(\ID). \]
Yet the condition $\Phi_h\in L^1_{loc}$ necessary to invoke Weyl's lemma (and of course the condition $\Phi\in L^1(\ID)$ is ubiqious in uniqueness theorems for extremal quasiconformal mappings after Hamilton's work \cite{Hamilton}) requires that at least 
 \[  \IK(w,h)\Psi'(\IK(w,h)) J(w,h) \lambda(h) \in L^1_{loc}(\ID) \]
This later condition is not going to occur automatically unless the growth of $\Psi$ is controlled.  Roughly we require
\[  t \Psi'(t) \approx \Psi(t) \]
This is fine in the case $\Psi(t) \approx t^p$ and so that is why there is always an Ahlfors-Hopf differential in that case. However, this growth is insufficient to guarantee topological regularity such as being locally homeomorphic or even openness.  To get this one needs $\Psi(t)\geq e^{\epsilon\,t}$.  However,  then the condition above cannot be satisfied.  One can now see where the hypotheses in Theorem \ref{main1} and the subsequent remark are coming from.  The situation is deftly avoided for mappings between closed Riemann surfaces by the Riemann-Roch theorem - the space of holomorphic Ahlfors-Hopf quadratic differentials is finite dimensional and an approximation technique can be used.  See \cite{MY3}.

\section{The inverse problem and Ahlfors-Hopf differentials.}

With the above discussion at hand,  our Theorem \ref{main1} reduces to the following results which are of independent interest.
The first results mentioned here are associated with the classical uniqueness properties of extremal quasiconformal mappings of Teichm\"uller type and are presumably well known to experts even in this more general setting.  However they lead us to the uniqueness theorems we seek and so include them here.

\begin{theorem}[mappings with harmonic argument]\label{thm1}
Let $S, R$ be closed Riemann surfaces. Let $f,g:S\to R$ be homotopic homeomorphisms of finite distortion, and suppose $\xi=g\circ f^{-1} \in W_{loc}^{1,1}(S,S)$. Suppose that $f$ has both the Lusin $\mathcal{N}$ and $\mathcal{N}^{-1}$ properties, and its Beltrami coefficient has the form 
\[
\mu_h(z)=|\mu_h(z)|\frac{\overline{\phi(z)}}{|\phi(z)|}, \quad h=f^{-1},
\]
where $\phi\in L^1(S)$ is a holomorphic quadratic differential.  If $|\mu_g|\leq|\mu_f|$, then $f=g$ in $S$.
\end{theorem}

Note this proves the uniqueness part of Teichm\"uller's theorem since if $f$ and $g$ are quasiconformal,  so is $g\circ f^{-1}$.

\begin{theorem}\label{thm2}
Let $(S,\sigma_S), (R,\sigma_R)$ be closed Riemann surfaces with their hyperbolic metrics. Let $h,H:S\to R$ be homotopic homeomorphisms of finite distortion with positive Jacobians, $J_h,J_H>0$ $a.e\; S$, and $\xi=H^{-1}\circ h\in W_{loc}^{1,1}(S)$. Let $\Psi:[1,\infty)\to[1,\infty)$ be a convex increasing function. Assume that $h$ has holomorphic Ahlfors-Hopf quadratic differential
\[
\Phi_h=\Psi'(\IK(z,h)) h_z\overline{h_\zbar} \, d\sigma_R(h) \in L^1(S),
\]
and
\[
\int_S \Psi(\IK(z,h)) \;d\sigma_R(h)   <\infty.
\]
Then
\begin{equation}\label{fs}
\int_S\Psi(\IK(z,h)) \;d\sigma_R(h)  \leq\int_S\Psi(\IK(z,H)) \;d\sigma_R(H)  
\end{equation}
Furthermore,  equality holds if and only if $h\equiv H$.
\end{theorem}

\medskip

\noindent{\bf Remark 2.} To make our notation clearer we make the following comments.  The continuous maps $h,H:S\to R$ lift to the universal covers as
$\hat{h},\hat{H}:\ID\to \ID$ and these maps are automorphic with respect to the Fuchsian groups $\Gamma_R=\pi_1(R)$ and $\Gamma_S=\pi_1(S)$.  The maps $\hat{h},\hat{H}$ induce the same isomorphism $\rho:\pi_1(S)\to \pi_1(R)$ between fundamental groups and 
\[ \hat{h}\circ \gamma = \rho(\gamma) \circ \hat{h}, \quad \hat{H} \circ \gamma = \rho(\gamma) \circ \hat{H}, \quad \quad \gamma\in \pi_1(S) \]
\[\begin{tikzcd}
\ID \arrow{r}{\hat{h}} \arrow[swap]{d}{\pi_S} & \ID \arrow{d}{\pi_R} \\
S \arrow{r}{h} & R
\end{tikzcd}
\]
The holomorphic differential $\phi$ lifts to a holomorphic function $\hat{\phi}:\ID \to \IC$
\[ \hat{\phi} = \Psi'(\IK(z,\hat{h})) \frac{\hat{h}_z\overline{\hat{h}_\zbar} }{(1-|\hat{h}|^2)^2} \]
and for {\em any} (say convex) fundamental polygon ${\cal P}$ for $\pi_1(S)$ we have
\[ \int_{\cal P}   \frac{|\hat{\phi}|dz}{(1-|z|^2)^2} <\infty \]
which in the context of {\em closed} surfaces is the same as $\hat{\phi}\in L^1_{loc}(\ID)$.  Similarly the ``finite energy" assumption 
$\int_S \Psi(\IK(z,h))   \;d\sigma_R(h)   <\infty $ should be read as
\[ \int_{\cal P} \Psi(\IK(z,\hat{h}))  \frac{J(z,\hat{h})dz}{(1-|\hat{h}|^2)^2}  <\infty \]
or equivalently $\Psi(\IK(z,\hat{h})) J(z,\hat{h} ) \in L^1_{loc}(\ID)$.  The final statement then reads as
\begin{equation}\label{fs2}
\int_{\cal P} \Psi(\IK(z,\hat{h}))\frac{J(z,\hat{h}) dz}{(1-|\hat{h}|^2)^2}  \leq\int_{\cal P} \Psi(\IK(z,\hat{H}))\frac{J(z,\hat{H}) dz}{(1-|\hat{H}|^2)^2}  
\end{equation}

We also consider the boundary problems on the planar disk. The same result applies and we get

\begin{theorem}\label{diskthm}
Let $f,g:\overline{\ID}\to\overline{\ID}$  both be finite distortion homeomorphisms, $f|_\IS=g|_\IS$, $H=g^{-1}$, $h=f^{-1}$, and $\xi=g\circ h\in W_{loc}^{1,1}(\ID)$. Assume that $f$ has both Lusin $\mathcal{N}$ and $\mathcal{N}^{-1}$ properties, and satisfies the equation
\[
\mu_h=|\mu_h|\frac{\overline{\phi}}{|\phi|},
\]
where $\phi\in L^1(\ID)$ is holomorphic, and $|\mu_g|\leq|\mu_f|$. Then $f=g$ in $\overline{\ID}$.  The requirement $\phi\in L^1(\ID)$ is necessary.
\end{theorem}

\begin{corollary}
Let $f:\overline{\ID}\to\overline{\ID}$ be a Teichm\"uller map, that is  its inverse $h=f^{-1}$ satisfies
\[
\mu_h=k\frac{\overline{\phi}}{|\phi|},\quad k\in [0,1),
\]
where $\phi\in L^1(\ID)$ is holomorphic. Let $g:\overline{\ID}\to\overline{\ID}$ be a quasiconformal mapping such that $g|_\IS=f|_\IS$, and $|\mu_g(z)|\leq k$ almost everywhere in $\ID$. Then $f\equiv g$.
\end{corollary}

\begin{theorem}\label{thmplanar}
Let $h,H:\overline{\ID}\to\overline{\ID}$ both be finite distortion homeomorphisms $h|_\IS=H|_\IS$, $J_h,J_H>0$ a.e., and $\xi=H^{-1}\circ h\in W_{loc}^{1,1}(\ID)$. Let $\Psi:[1,\infty)\to[1,\infty)$ be a convex function. Assume that $h$ has holomorphic Hopf differential
\[
\Phi_h=\Psi'(\IK(z,h))\,h_z\overline{h_\zbar} \in L^1(\ID),
\]
and
\[
\int_\ID\Psi(\IK(z,h))\,J(z,h) \; dz <\infty.
\]
Then
\[
\int_\ID\Psi(\IK(z,h))\,J(z,h)\; dz \leq\int_\ID\Psi(\IK(z,H))\,J(z,H)\; dz
\]
Furthermore,  equality holds if and only if $h\equiv H$.
\end{theorem}

\section{The Reich-Strebel inequalities}
\begin{lemma} \label{RS} (Reich-Strebel inequality) Let $S$ be a compact Riemann surface. Let $\varphi\in L^1(S)$ be a holomorphic quadratic differential on $S$. Let $f$ be a finite distortion self-homeomorphism of $S$ which is homotopic to the identity. Then 
\begin{equation}\label{1}
\int_S|\phi|\leq\int_S\sqrt{|\phi(f)|}\sqrt{|\phi|}\left|f_z-\frac{\phi}{|\phi|}f_\zbar\right|.
\end{equation}
\end{lemma}

We remark the classic Reich-Strebel inequality was stated with the assumption that $f$ is a quasiconformal mapping, see \cite{RS2}, or \cite{G1}. However, the same argument applies if $f$ is only a mapping of finite distortion.

Note that (\ref{1}) leads to the following two inequalities by use of the  Cauchy-Schwarz inequality:
\begin{equation}\label{2}
\int_S|\phi(f)|\left|f_z-\frac{\phi}{|\phi|}f_\zbar\right|^2\geq\frac{\int_S\left(\sqrt{|\phi(f)|}\sqrt{|\phi|}\left|f_z-\frac{\phi}{|\phi|}f_\zbar\right|\right)^2}{\int_S|\phi|}\geq\int_S|\phi|,
\end{equation} 
\begin{equation}\label{3}
\int_S|\phi|\frac{\left|f_z-\frac{\phi}{|\phi|}f_\zbar\right|^2}{J(z,f)}\geq\frac{\int_S\left(\sqrt{|\phi(f)|}\sqrt{|\phi|}\left|f_z-\frac{\phi}{|\phi|}f_\zbar\right|\right)^2}{\int_S|\phi(f)|J(z,f)}\geq\int_S|\phi|.
\end{equation}
\bigskip

There is also a version of this inequality for finite distortion mappings on planar domains, see \cite{IO1}.

\begin{lemma}[Iwaniec-Onninen]
Let $f:\overline{\ID}\to\overline{\ID}$ be a finite distortion homeomorphism with $f|_\IS=Identity|_\IS$ and let $\phi\in L^1(\ID)$ be holomorphic. Then
\begin{equation}\label{4}
\int_\ID|\phi|\leq\int_\ID\sqrt{|\phi(f)|}\sqrt{|\phi|}\left|f_z-\frac{\phi}{|\phi|}f_\zbar\right|.
\end{equation}
\end{lemma}

\bigskip

\section{Proofs}
We now prove Theorem \ref{thm1} and Theorem \ref{thm2}. We illustrate why the proofs for the other results are entirely similar.\\

\noindent{\bf Proof of Theorem 1.} From (\ref{2}) we have
\begin{equation}\label{5}
\int_S|\phi(\xi)|\left(\left|\xi_w-\frac{\phi}{|\phi|}\xi_\wbar\right|^2-J(z,\xi)\right)\geq0.
\end{equation}
However, we also have the following pointwise inequality:
\begin{equation}\label{6}
\left|\xi_w-\frac{\phi}{|\phi|}\xi_\wbar\right|^2\leq J(z,\xi).
\end{equation}
To see this, we can rewrite this inequality in the following way.
\begin{align*}
0&\geq\left|\xi_w(f)-\frac{\phi(f)}{|\phi|(f)}\xi_\wbar(f)\right|^2-J(f,\xi)\\
&=\left|(g\circ h)_w(f)-\frac{\phi(f)}{|\phi(f)|}(g\circ h)_\wbar(f)\right|^2-J(f,{g\circ h})\\
&=\left|(1-|\mu_f|)\left(\frac{g_z\overline{f_z}}{J(z,f)}-\frac{\phi(f)}{|\phi(f)|}\frac{g_\zbar f_z}{J(z,f)}\right)\right|^2-\frac{J(z,g)}{J(z,f)}.
\end{align*}
Also,
\begin{align}
\left|(1-|\mu_f|)\left(\frac{g_z\overline{f_z}}{J(z,f)}-\frac{\phi(f)}{|\phi(f)|}\frac{g_\zbar f_z}{J(z,f)}\right)\right|^2\leq&\frac{(1-|\mu_f|)^2}{J(z,f)^2}(|g_z|+|g_\zbar|)^2|f_z|^2\notag\\
=&\frac{1-|\mu_f|}{1+|\mu_f|}\frac{1+|\mu_g|}{1-|\mu_g|}\frac{J(z,g)}{J(z,f)}\leq\frac{J(z,g)}{J(z,f)}\label{7}.
\end{align}
This proves (\ref{6}), and then together with (\ref{5}) we have
\[
\left|(1-|\mu_f|)\left(\frac{g_z\overline{f_z}}{J(z,f)}-\frac{\phi(f)}{|\phi(f)|}\frac{g_\zbar f_z}{J(z,f)}\right)\right|^2=\frac{J(z,g)}{J(z,f)}.
\]
Next,  in the computation of (\ref{7}) we see that equality holds only when $\mu_f=\mu_g$. Together with the condition that $f$ and $g$ are homotopic we find that $f=g$ in $S$.\hfill$\Box$

\bigskip

\noindent{\bf Proof of Theorem 2.} We consider $H=h\circ\xi^{-1}$. The composition formula for the distortion function $\IK_H$ follwows from the usual composition formula for the Beltrami coefficients of quasiconformal mappings and gives us the following lemma by direct calculation.
\begin{lemma} With the notation above,
\begin{equation}
\IK(\xi,H)=\IK(z,\xi)\IK(z,h)\left[1-\frac{4\Re e(\mu_\xi\overline{\mu_h})}{(1+|\mu_\xi|^2)(1+|\mu_h|^2)}\right].
\end{equation}
\end{lemma} 

We make the following calculation on the surface $S$ by considering the universal cover so the reader can plainly see why the argument also works for the disk.  Thus we set
\[ \lambda(z) = \frac{1}{(1-|z|^2)^2} \]
the hyperbolic area element in the disk.  We will not use the particular form of $\lambda$ in what follows, and this is why the calculation also holds in the disk.  We calculate with ${\cal P}$ a  convex fundamental polyhedron and write $\xi$ for $\hat{H}^{-1}\circ \hat{h}$.
\begin{eqnarray*}
\lefteqn{\int_S\Psi(\IK(z,\hat{H})) d\sigma_R(\hat{H}) }\\ &= & \int_{\cal P} \Psi(\IK(z,\hat{H}))\, \lambda(\hat{H}) \; dz = \int_{\cal P}\Psi(\IK(\xi,\hat{H})) J(\xi,\hat{H})J(z,\xi) \lambda(\hat{H}) \\
&= & \int_{\cal P} \Psi\left(\IK(z,\xi)\IK(z,\hat{h}) \left[1-\frac{4\Re e(\mu_\xi \overline{\mu_{\hat{h}}})}{(1+|\mu_\xi|^2)(1+|\mu_{\hat{h}}|^2)}\right]\right) \, J(z,\hat{h})\,\lambda(\hat{h}).
\end{eqnarray*}
Then convexity of $\Psi$ gives us that
\[
\Psi(y)-\Psi(x)\geq(y-x)\Psi'(x),
\]
and we can make the following calculation.
\begin{eqnarray*}
\lefteqn{\int_{\cal P} \Psi(\IK(z,\hat{H}))\, J(z,\hat{H})\, \lambda(\hat{H}) -\Psi(\IK(z,\hat{h}))\,J(z,\hat{h})\, \lambda(\hat{h}) } \\ 
&= &\int_{\cal P} \left[\Psi\left(\IK(z,\xi)\IK(z,\hat{H})\left[1-\frac{4\Re e(\mu_\xi\overline{\mu_{\hat{h}}})}{(1+|\mu_\xi|^2)(1+|\mu_{\hat{h}}|^2)}\right]\right)-\Psi(\IK(z,\hat{h}))\right]\,J(z,\hat{h})\, \lambda(\hat{h}) \\
& \geq&\int_{\cal P} \left(\IK(z,\xi) \left[1-\frac{4\Re e(\mu_\xi\overline{\mu_{\hat{h}}})}{(1+|\mu_\xi|^2)(1+|\mu_{\hat{h}}|^2)}\right]-1\right)\IK(z,\hat{h})\Psi'(\IK(z,\hat{h}))\, J(z,\hat{h}) \,\lambda(\hat{h})\\
&=&2\int_{\cal P} \frac{(|\hat{h}_z|-|\hat{h}_\zbar|)^2|\xi_\zbar|^2}{J(z,\xi)} \Psi'(\IK(z,\hat{h}))\lambda(\hat{h}) \\ && 
\quad +2\int_{\cal P} \left(\frac{|\xi_z-\frac{\phi}{|\phi|}\xi_\zbar|^2}{J(z,\xi)}-1\right)|\phi|\geq 0.
\end{eqnarray*}
Here the last inequality follows from (\ref{3}). Furthermore, if  equality holds, then both of the terms in the last equality must be $0$. In particular, the first term
\[
\int_{\cal P} \frac{(|\hat{h}_z|-|\hat{h}_\zbar|)^2|\xi_\zbar|^2}{J(z,\xi)}\,\Psi'(\IK(z,\hat{h}))\lambda(\hat{h}) =0.
\]
Since $|\hat{h}_z|>|\hat{h}_\zbar|$ almost everywhere,  this holds only if $\xi_\zbar=0$. That is $\xi$ is conformal and since it is homotopic to the identity, $\xi\equiv Identity$ (in $S$). Then we see that  $\hat{h}\equiv \hat{H}$. \hfill$\Box$
 
\section{Connections.}

Now if $f$ is a minimiser (or even a variational critical point) to problems A or B and satisfies the hypothesis of Theorem \ref{main1},  then $f^{-1}=h$ is a diffeomorphic variational critical point for the functional
\[ h\mapsto \int \Psi'(\IK(z,h))\, J(z,h)\, \lambda(h) \]
and as such admits a holomorphic Ahlfors-Hopf differential. The previous sections shows $h$ to be unique, and thus $f$ is unique.

\section{Mappings with boundary values equal to the identity.} 

In this section we outline examples to show that the hypotheses of Theorem \ref{thmplanar} are necessary.  It is quite difficult to explicitly construct solutions to the Ahlfors-Hopf equation on the disk $\ID$ and so we move the problem to the upper-half space $\IH^2$.   Define $\psi:\ID \to \IH^2$ by
\[ \psi(z) = -i\, \frac{z-1}{z+1}, \quad 1\mapsto 0, \; -1\mapsto \infty,  \; 0\mapsto i. \]
\begin{lemma} Suppose for a positive weight $\eta$ the mapping $h:\IH^2\to\IH^2$, $h\in W^{1,2}_{Loc}(\IH^2)$, has holomorphic Ahlfors-Hopf differential
\begin{equation}
\Phi_h=\Psi(\IK(z,h)) h_z \overline{h_\zbar} \; \eta(h) .
\end{equation}
Set
\begin{equation}
 \lambda(w)=\eta(\psi(w)) |\psi'(w)|^{2}
\end{equation}
Then $g=\psi^{-1}\circ h \circ \psi:\ID \to \ID$ has holomorphic Ahlfors-Hof differential
\begin{equation}
\Psi'(\IK(w,g))g_w\overline{g_\wbar}\;\lambda(g) = \Phi(\psi)(\psi')^2
\end{equation}
\end{lemma}
\noindent{\bf Proof.} We simply calculate
\begin{eqnarray*}
\lefteqn{\Psi'(\IK(w,g))g_w\overline{g_\wbar}\;\lambda(g)}\\
&=& \Psi'(\IK(w,h\circ\psi)) |(\psi^{-1})'(h\circ \psi)|^2 h_z(\psi) \overline{h_\zbar(\psi)} (\psi')^2 \;\lambda(g)\\
&=& \Phi(\psi)  (\psi')^2   \; |(\psi^{-1})'(h\circ \psi)|^2 \;\lambda(\psi^{-1}\circ h\circ \psi)/\eta(h\circ\psi)
\end{eqnarray*}
and this last term is equal to $1$. \hfill $\Box$

\medskip

\begin{lemma}Should $\eta(z)=\Im m(z)^{-2}$,  the hyperbolic area element of $\IH^2$,  then $\lambda(w)=(1-|w|^2)^{-2}$,  the hyperbolic area element of the disk.
\end{lemma}

The first examples to consider are the linear mappings $h(x,y)=x+i \;\alpha \,y$ with $\alpha>0$ and $\eta\equiv 1$.  Then 
\[ \Psi(\IK(z,h)) h_z \overline{h_\zbar} \; \eta(h)  = \frac{1}{4} \Psi\big(\frac{\alpha+1/\alpha}{2}\big)(1-\alpha^2) = \beta \]
is constant,  and so holomorphic.  Then $g_\alpha=\psi\circ h \circ \psi^{-1}:\ID\to\ID$ is a diffeomorphism with identity boundary values and holomorphic Ahlfors-Hopf differential with respect to the weight 
\[ \lambda(z)= |\psi'(w)|^{2} = \frac{4}{|z+1|^4} \] 
and equal to 
\[ \Phi_{g_\alpha}= \Psi(\IK(w,g_\alpha)) (g_\alpha)_w \overline{(g_\alpha)_\wbar} \; \lambda(g_\alpha)  = \beta\;  (\psi')^2 \] 
unless $\alpha=1$ whereupon $g_\alpha$ is the identity and $\Phi_{g_\alpha}=0$.

\medskip
All these mappings $g_\alpha$ are quasiconformal,  the identity on the boundary of $\ID$ and have holomorphic Ahlfors-Hopf differential $\Phi_{g_\alpha}$ with respect to the same weight.  Only one of them,  $g_1$,  has $\Phi_{g_\alpha}\in L^1(\ID)$ however.  This shows the hypothesis $\Phi_{g}\in L^1(\ID)$ in Theorem \ref{thm2} to be necessary.

\medskip

We next consider the case of the two hyperbolic metrics. In fact we shall assume $\eta(z)$ depends only on the imaginary part of $z$ and defines a complete metric.  In particular $\eta(s)\to \infty$ as $s\to 0$.  Following the above example we set
\begin{equation}
h(z)=h(x+iy)=x+i u(y),  \quad u(0)=0,\; u'(0)>0, \lim_{a\to\infty} u(a)=+\infty.
\end{equation}
We calculate
\[ h_z = \frac{1}{2}(1+u'(y)),\quad h_\zbar = \frac{1}{2}(1-u'(y)), \quad \IK(z,h)=\frac{1+u'(y)^2}{2u'(y)}.\]
We seek $h$ with constant Ahlfors-Hopf differential equal to $\lambda$.  This yields the following ordinary differential equation for $u(y)$.
\begin{equation}\label{ode}
\Psi'\Big(\frac{1+u'(y)^2}{2u'(y)}\Big)(1-u'(y)^2) \, \eta(u(y)) = 4\lambda.
\end{equation}
Since $\lambda$ is assumed constant there are two classes of solutions.  Those with $0<u'<1$ and those with $1<u'<\infty$ characterised by $\lambda>0$ and $\lambda<0$ respectively. Of course $\lambda=0$ gives the identity. Consider the function 
\begin{equation}\label{convex} t\stackrel{{\cal F}}\mapsto \Psi'\Big(\frac{1+t^2}{2t}\Big)(1-t^2), \quad  {\cal F}(1)=0.  \end{equation}
The derivative of ${\cal F}$ is
\[  {\cal F}'(t)=\frac{1}{2}\Psi''\Big(\frac{1+t^2}{2t}\Big)(1-t^2)(1-\frac{1}{t^2}) -2t\Psi'\Big(\frac{1+t^2}{2t}\Big)\leq 0, \quad t\neq 0\]
Since $\Psi''\geq 0$ and $\Psi'>0$, ${\cal F}'(t)=0$ implies  $t=0$ and so the function defined by (\ref{convex}) is strictly decreasing.
\begin{itemize}
\item $\lambda>0$, $t\in (0,1)$. In this case ${\cal F}$ decreases from 
\begin{equation}\label{Mdef}
M=\lim_{a\to 0}{\cal F}(a)
\end{equation}
 to $0$. Often, but not always $M=+\infty$.
\item $\lambda<0$, $t\in (1,\infty)$. In this case ${\cal F}$ decreases from $0$ to $-\infty$.
\end{itemize}
This last assertion follows since we cannot have $\Psi'\Big(\frac{1+t^2}{2t}\Big)t^2$ bounded for convex increasing $\Psi$.
We also remark that as $\eta\neq 0$ the continuous function ${\cal F}$ has constant sign if we are to satisfy (\ref{ode}).

\medskip

 If $\big\{ {\lambda}/{\eta(u)}:u\in (0,\infty)\big\}$, lies in the range of ${\cal F}$ we can rewrite (\ref{ode}) as
\begin{equation}\label{ode2}
u' = {\cal G}(u)={\cal F}^{-1}\Big(\frac{\lambda}{\eta(u)}\Big) >0,\quad u(0)=0.
\end{equation}
where ${\cal G}$ is a smooth function of  $u$.  Notice we already see that if $\lambda>0$, $M<\infty$ (defined at (\ref{Mdef})) and if $\eta(u)\to 0$ as $u\to+\infty$,  then there cannot be a solution to the o.d.e.  This happens in the hyperbolic case $\eta(y)=y^{-2}$. The separable equation (\ref{ode2}) admits the solution
\[ u(y)=\int_{0}^{y} {\cal G}(t) dt \]
There is only one further requirement we must check and that is that
$ \int_{0}^{\infty} {\cal G}(t)\; dt =+\infty $
ensuring $h$ is surjective.  Notice that ${\cal F}^{-1}(0)=1$ and so there is never a singularity at $0$. We have  
\begin{eqnarray}\label{divergence}
\int_{0}^{x} {\cal G}(t)\; dt  &=& \int_{0}^{x} {\cal F}^{-1}\Big(\frac{\lambda}{\eta(t)}\Big) \; dt  
\end{eqnarray}
We consider  two cases.

\subsubsection{$\lambda > 0$} In this case we must have $0<u'<1$ and ${\cal F}\geq 0$. The range of ${\cal F}$ is the interval $(0,M)$. We leave the reader to consider the case $M<\infty$ as an interesting technical diversion. With $M=\infty$ we can always write down a solution and have ${\cal F}(0)=\infty$. If 
\[ {\cal F}^{-1}\Big(\frac{\lambda}{\eta(s)}\Big) \geq \frac{\delta}{s}, \quad \delta>0, \;\; s>>1, \]
then the divergence of the integral at (\ref{divergence}) as $x\to\infty$ is assured. As ${\cal F}$ is strictly decreasing so too is ${\cal F}^{-1}$ and so the estimate above is the same as the estimate
\[ \frac{\lambda}{\eta(s)} \leq {\cal F}\Big(\frac{\delta}{s} \Big)= \Psi'\Big(\frac{1}{2}\big(\frac{s}{\delta}+\frac{\delta}{s} \big)\Big)(1-\frac{\delta^2}{s^2}),\quad \delta>0, \; s>>1.\]
If the quantity $\eta(s) \Psi'(s)$ is bounded below as $s\to \infty$,  then since $\Psi$ is convex we can chose a positive $\delta$ so that we achieve divergence.  

\begin{lemma}  Suppose that $\lambda>0$, that $\lim_{t\to 0} \Psi'\Big(\frac{1+t^2}{2t}\Big)=+\infty$ and that 
\[ \liminf_{s\to\infty} \eta(s) \Psi'(s) \geq \epsilon >0.\] Then there is a diffeomorphism $h:\IH^2\to\IH^2$,  $h|\partial \IH^2=identity$ and
\[ \Psi'(\IK(z,h)) h_z \overline{h_\zbar}\; \eta(\Im m(h))  =\lambda.\]
If $\eta(s)$ is bounded below these solutions are quasiconformal.
\end{lemma}
For the hyperbolic metric,  the growth condition is $\Psi'(s)\geq c s^2$.  With the canonical examples $\Psi(s)=s^p$ ($p>1$ so $M=+\infty$) we see that there is always a solution if $p\geq 3$ in this case.  However,  if $p=2$,
\[ {\cal F}(t)= \frac{1-t^4}{t},\quad {\cal F}^{-1}(\lambda t^{2}) \approx \frac{1}{\lambda  t^{2}}, \]
and the integral converges.  Thus there is no solution $h$ of the form we propose,  and probably no solution at all.  Further notice that  $\eta(s)$ is bounded below if and only if $u'$ is bounded below and $h$ is quasiconformal.  This can never be the case for the hyperbolic metric.

\subsubsection{$\lambda < 0$} In  this case we have $u'\geq 1$ and the range of ${\cal F}$ is $(-\infty,0)$.  We can always solve the equation (\ref{ode}) and there is no question of divergence since ${\cal F}$ decreases to $-\infty$. 

\begin{lemma}  Suppose that $\lambda<0$.  Then there is a diffeomorphism $h:\IH^2\to\IH^2$,  $h|\partial \IH^2=identity$ and
\[ \Psi'(\IK(z,h)) h_z \overline{h_\zbar}\; \eta(\Im m(h))  =\lambda.\]
If $\eta(s)$ is bounded below these solutions are quasiconformal.
\end{lemma}

\subsection{Hyperbolic Harmonic case.} We briefly consider what happens in the case of harmonic mappings in the hyperbolic metric,  where $\eta(s)=s^{-2}$ and also $\Psi(t)=t$.  Then we have
\begin{equation}\label{odeharmonic}
 \sqrt{1-4 \lambda u(y)^2}= u'(y), \quad u(0)=0, \; u'(0)=1.
\end{equation}
The solution for $\lambda<0$ is
\[ u(s) = \int_{0}^s \frac{du}{\sqrt{1-4 \lambda u^2}} = \frac{1}{2\sqrt{|\lambda|}} {\rm arcsinh} (2\sqrt{|\lambda|} s) \]
The distortion of the induced mapping $h$ is 
\[ \IK(w,h)= \frac{2 \lambda  s^2+1}{\sqrt{4 \lambda  s^2+1}}, \quad s=\Im m(w) \]
In particular this mapping is not quasiconformal unless $\lambda=0$ and $h=identity$.

\subsection{$\Psi(t)=t^2$} In this case we have
\begin{equation}\label{ode}
 1-u'^4 = 4  \lambda u' u^2,  \quad u(0)=0, \; u'(0)=1.
\end{equation}
If $\lambda >0$ the left-hand side is bounded above by $1$ and the right-hand side is positive we must have $u'\to 0$ and $u' u^2\approx c$ so that $u(s)\approx c s^{1/3}$.  While if $\lambda<0$ we have $u'^4 \approx 4  \lambda u' u^2$ and $u \approx  c s^3$.

\end{document}